\documentclass[12pt]{article}
\usepackage{amsmath}
\usepackage{amsthm}
\usepackage{amsfonts}
\usepackage{amssymb}

\newtheorem{theorem}{Theorem}
\newtheorem{proposition}[theorem]{Proposition}

\newtheorem{corollary}[theorem]{Corollary}

\theoremstyle{definition}
\newtheorem{definition}[theorem]{Definition}
\newtheorem{example}[theorem]{Example}

\theoremstyle{remark}
\newtheorem{remark}[theorem]{\bf Remark}
\newcommand {\Z} {{\mathbb  Z}}
\newcommand {\R} {{\mathbb  R}}
\newcommand {\N} {{\mathbb  N}}
\allowdisplaybreaks[4]
\begin{document}
\title{Large classes of minimally supported frequency wavelets of 
$L^2(\R)$ and $H^2(\R)$}
\author{
Nicola Arcozzi
\thanks{Department of Mathematics, 
University of Bologna, 
Piazza di Porta San donata, 5; 
40127, ITALY. 
Email: arcozzi@dm.unibo.it}
\and 
Biswaranjan Behera
\thanks{Department of Mathematics,
Indian Institute of Technology, 
Kanpur 208016, INDIA.
Current address: 
Stat-Math Unit, 
Indian Statistical Institute, 
203 B. T. Road, Calcutta, 700108, 
India. 
Email: br\_behera@yahoo.com
Research supported by the National Board 
for Higher Mathematics, Govt.\ of India.}
\and 
Shobha Madan
\thanks{Department of Mathematics, 
Indian Institute of Technology, 
Kanpur 208016, INDIA. 
Email: madan@iitk.ac.in}
}
%\author{Nicola Arcozzi\\Dipartimento di Matematica\\Universita di %Bologna\\ITALY
%\and Biswaranjan Behera\\Department of Mathematics\\Indian Institute of Technology\\Kanpur,INDIA
%\and Shobha Madan\\Department of Mathematics\\Indian Institute of Technology\\Kanpur, INDIA}
\maketitle
%\tableofcontents

\begin{abstract}
We introduce a method to construct large classes of MSF wavelets
of the Hardy space $H^2(\R)$ and symmetric MSF wavelets of $L^2(\R)$, 
and discuss the classification of such sets. As application, we show 
that there are uncountably many wavelet sets of $L^2(\R)$ and 
$H^2(\R)$. We also enumerate all symmetric wavelets of $L^2(\R)$ with 
at most three intervals in the positive axis as well as 3-interval 
wavelet sets of $H^2(\R)$. Finally, we construct families of MSF
wavelets of $L^2(\R)$ whose Fourier transform does not vanish in any
neighbourhood of the origin.
\end{abstract}

\vskip 3mm
\begin{footnotesize}
\noindent{\bf Notation.} 
In this article, {\it measure} will always mean Lebesgue measure. 
All subsets of the real line we consider,
are assumed to be Lebesgue measurable. With $|E|$ we denote the
Lebesgue measure of a measurable set $E$ in $\mathbb R$. We say
that a relation between measurable sets {\it holds almost
everywhere} ($a.e.$) if their characteristic functions are equal
$a.e$. Thus, $A=B\ a.e.$ means that $\chi_{_A}=\chi_{_B}\ a.e.$,
and $\coprod_{n}A_n=A\ a.e.$ means that
$\sum_n\chi_{_{A_n}}=\chi_{_A}\ a.e.$, where $\coprod$ denotes the
disjoint union.
\end{footnotesize}
%%================================================================

\section{Introduction}
This article is an attempt to better understanding MSF (i.e.,
Minimally Supported Frequency) wavelets of $L^2(\R)$ and
$H^2(\R)$. Our results have short and elementary proofs. They
include geometric algorithms to construct MSF wavelets, and
several new examples of MSF wavelets, a class of which answers in
the negative a question of Eugenio Hern\'andez. Another family of
wavelets constructed provides a positive answer to a question of
G. Garrig\'os regarding the existence of wavelets of $L^2(\R)$
with certain properties.

Recall that an {\it orthonormal wavelet} is a function
 $\psi\in L^2({\mathbb R})$ such that
 $\{\psi_{j,k}\colon j,k\in{\mathbb Z}\}$ is an orthonormal
 basis for $L^2(\mathbb R)$, where
\[
 \psi_{j,k}(x)=2^{j/2}\psi(2^jx-k).
\]
The wavelet $\psi$ is called an {\it MSF wavelet} if there exists 
a set $K$ in ${\mathbb R}$ such that
\begin{equation} \label{eqpsi}
 \left|\widehat{\psi}\right|=\chi_{_K}.
\end{equation}
In this case, we say that $K$ is a {\it wavelet set}. As
definition of Fourier transform, we use
\[
\widehat{\psi}(\xi)=
\int_{-\infty}^{+\infty}\psi(x)\exp({-2\pi ix\xi})dx
\]
which is different from that used, for instance, in~\cite{HW},
but which highlights the essentially number-theoretic nature of
our considerations. To compare our sets with the ones in~\cite{HW}, 
it suffices to multiply our subsets on the real line by
a factor $2\pi$.

The motivation for the name MSF is that orthonormal wavelets have
the property that $|\mbox{supp}~\hat{\psi}|\ge 1$, with equality
if and only if~(\ref{eqpsi}) holds $a.e.$ for some measurable $K$
in $\mathbb R$. If $K$ is a wavelet set, all functions $\psi$ such
that
\[
 \left|\widehat{\psi}(\xi)\right|=\chi_K(\xi)\ a.e.
\]
are orthonormal wavelets. For these facts and the theorem below,
see~\cite{HW}, especially \S7.2.

The wavelet sets admit a simple characterization in terms of two
geometric conditions, that will be our starting point.

\begin{theorem}[\cite{HW}, \S7.2, Theorem 2.3]
\label{Tmsf}
Let $K$ be a subset of ${\mathbb R}$. $K$ is a wavelet set if and only if the
following are satisfied
\begin{enumerate}
\item[{\rm (T)}] $\coprod_{n\in {\mathbb Z}}(K+n)={\mathbb R}\ a.e.$
\item[{\rm (D)}] $\coprod_{n\in{\mathbb Z}}2^nK={\mathbb R}\ a.e.$
\end{enumerate}
\end{theorem}

Conditions (T) and (D) of Theorem~\ref{Tmsf} specialize to the MSF
case the four conditions that, in full generality, characterize
the Fourier transform of an orthonormal wavelet. See~\cite{HW},
Chapter $7$, for a thorough discussion of this issue. As in the
general theory, in (T) and (D) the additive (translation by
integers) and the multiplicative (dilation by powers of 2)
structure of the real line are competing.

Wavelet sets are plentiful and are perhaps beyond the possibility
of a simple and complete classification. A natural subfamily of
them is the class of {\it symmetric wavelet sets}, i.e., sets $K$
such that
\begin{equation}
\label{eqsym}
K\cap (-\infty,0]=-(K\cap[0,\infty)).
\end{equation}

From now on, we denote by $K^-$ and $K^+$, respectively, the sets
appearing on the left and right hand side of~(\ref{eqsym}).
Symmetric wavelet sets are completely classified when $K^+$ is the
union of finitely many disjoint intervals
\begin{equation}
\label{eqKplus}
K^+=I_1\cup\dots\cup I_n
\end{equation}
and $n=1,2$~(\cite{HKLS}). 
%Some examples exist in the literature. 
In this note, we present a method to construct large families of such
sets for all $n\ge 1$.

One construction consists in associating symmetric wavelet sets to certain 
polygons having vertices in a dyadic lattice. In this way, we obtain, 
for each $n$, a family of wavelet sets that depends on a finite 
number of arbitrary parameters with values in the positive integers.
For $n=1,2$, this provides a complete characterization of the 
symmetric wavelet sets. For $n=1$, the only symmetric wavelet set is the 
{\it Shannon set}, for which $K^+=\left[\frac{1}{2},1\right]$.
%\[
%K^+=\left[\frac{1}{2},1\right].
%\]
In the case $n=2$, the symmetric wavelet sets form a family depending on
one integer parameter, which was first identified in~\cite{HKLS}.

In the case $n=3$, the wavelet sets depend on several parameters,
and they exhibit behaviours so far unnoticed. For instance, there
exists a family of symmetric wavelet sets in which the endpoints
of the intervals continuously depend on a real parameter. Hence,
for $n=3$, there are uncountably many symmetric wavelet sets. The
sets in this family are those for which
\begin{equation}
\label{eqEugenio}
K_a^+=\left[a,\tfrac{1}{2}\right]\cup\left[1-a,2a\right]
\cup\left[1,2(1-a)\right]
\end{equation}
where $a\in(\tfrac{1}{3},\tfrac{1}{2})$.
In particular, there exist symmetric wavelet sets in which more than four
endpoints {\it do not have} the form 
$\frac{p}{2^q},\ p\in{\mathbb N},\ q\in{\mathbb Z}$. 
This answers in the negative a question of Eugenio Hern\'andez.
%%=============================================================
\subsection*{Wavelets for the Hardy space $H^2(\R)$}

The classical Hardy space $H^2(\R)$ is the collection of all
functions of $L^2(\R)$ whose Fourier transform is supported in
$\R^+ = (0,\infty)$:
\[
H^2(\R) = \{f \in L^2(\R) : \hat f (\xi)=0 ~{\rm for~a.e.}~\xi
\leq 0 \}.
\]

It is clear that $H^2(\R)$ is a closed subspace of $L^2(\R)$. As
in the case of $L^2(\R)$, we can define a wavelet for $H^2(\R)$.
\textit{A function $\psi\in H^2(\R)$ is said to be a wavelet of
$H^2(\R)$ if the system of functions
\[
\{\psi_{j,k}:=2^{j/2}\psi(2^j\cdot-k): j, k\in \Z\}
\]
forms an orthonormal basis for $H^2(\R)$.} Such a function $\psi$
will be called an $H^2$-wavelet. An example of an $H^2$-wavelet is
the function whose Fourier transform is the characteristic
function of the interval $[1,2]$. In fact, for a long time this
was the only known \mbox{$H^2$-wavelet}. P. Auscher~\cite{aus}
proved that if $\psi\in H^2(\R)$ is such that $|\hat\psi|$ is
continuous on $\R$, and 
$|\hat\psi (\xi)|=O((1+|\xi|)^{-\alpha-\frac{1}{2}})$ at $\infty$, 
for some $\alpha>0$, then $\psi$ cannot be an $H^2$-wavelet. 
In particular, there is no band-limited $H^2$-wavelet such that 
$|\hat\psi|$ is continuous (A function $f$ is said to be band-limited 
if $\hat f$ has compact support).

An $H^2$-wavelet $\psi$ will be called an {\it $H^2$-MSF wavelet}
if $|\hat\psi|=\chi_K$ for some measurable subset $K$ of $\R^+$.
The associated set $K$, which has measure $1$, will be called an
{\it $H^2$-wavelet set}. The following theorem characterizes all
$H^2$-wavelet sets.

\begin{theorem}
\label{thm:h2wav} 
A set $K\subset \R^+$ is an $H^2$-wavelet set if and only if
\begin{itemize}
\item[{\rm (T')}] $\coprod_{k\in\Z}(K+k)=\R$ a.e.
\item[{\rm (D')}] $\coprod_{j\in\Z}2^j K=\R^+$ a.e.
\end{itemize}
\end{theorem}

The proof of Theorem \ref{thm:h2wav} can be obtained from the
corresponding theorem (Theorem~\ref{Tmsf}) in the usual $L^2(\R)$
case, with necessary modifications.

We shall call an $H^2$-wavelet set an {\it interval $H^2$-wavelet set} 
if it is the union of a finite number of
intervals of $\R^+$. In~\cite{HKLS} the authors characterized all
$H^2$-wavelet sets consisting of at most two intervals. In fact,
the only $H^2$-wavelet set which is a single interval is $[1,2]$,
and those which are union of two disjoint intervals are the
following sets:
\[
 \left[\frac{k+1}{2^{r+1}-1},\frac{k}{2^r-1}\right]
 \cup\left[\frac{2^rk}{2^r-1},\frac{2^{r+1}(k+1)}{2^{r+1}-1}\right],
\]
 where $r>0$, $0<k<2(2^r-1)$; $r,k\in\Z$.

\medskip

This article is organized as follows. In section~2 we present a
geometric construction for some symmetric wavelet sets of
$L^2(\R)$. In section~3 we discuss the problem of the complete
classification of the symmetric wavelet sets of $L^2(\R)$, with
some preliminary results. The geometric construction of section~2 is
extended to the case of $H^2(\R)$ in section~4. In section~5 we prove 
a result on the structure of interval $H^2$-wavelet sets. In section~6 
we give some examples and applications. In particular, we characterize
3-interval $H^2$-wavelet sets and symmetric wavelet sets $K$ of
$L^2(\R)$ for which $K^+$ is a union of three disjoint intervals.
We also show that there are uncountably many wavelet sets of $H^2(\R)$ 
and symmetric wavelet 
sets of $L^2(\R)$. In the last section we construct
three families of wavelet sets of $L^2(\R)$ having the origin as
an accumulation point. These constructions show the existence of
such wavelet sets which are (a) symmetric and unbounded 
(b) bounded and non-symmetric, and (c) bounded and symmetric.

Finally, let us note that the method presented here can be used to
construct non-symmetric wavelet sets of $L^2(\R)$, as well. In
this case, however, the construction only covers a very small
family of the wavelet sets.
%%================================================================
\section{A geometric construction of symmetric\\wavelet sets of $L^2(\R)$}

In this section, we provide a geometric algorithm to construct
symmetric wavelet sets of $L^2(\R)$, which is summarized in Theorem~\ref{Tpol}.

In the first quadrant of the Cartesian plane, consider the set $D$
of the points $P$ such that, for some
 $m\in{\mathbb N}_0={\mathbb N}\cup\{0\},\ \lambda\in {\mathbb Z}$,
\[
 P \equiv P[\lambda,m]=(2^{-\lambda},2^{-\lambda}m).
\]
Let ${\mathcal P}=(P_1\dots,P_n)$ be an ordered sequence of points
in $D$, $P_j=P[\lambda_j,m_j]$. For $j=1,\dots,n-1$, let
\begin{equation}
\label{eqaj}
a_j=-\frac{m_j2^{-\lambda_j}-m_{j+1}2^{-\lambda_{j+1}}}
{2^{-\lambda_j}-2^{-\lambda_{j+1}}},
\end{equation}
i.e., the negative of the slope of the straight line through $P_j$ and
$P_{j+1}$. We say that ${\cal P}$ is an {\it MSF polygonal} if
\begin{equation}
\label{eqMSFo} \lambda_1=0, \quad
4m_1=2^{-\lambda_n}(2m_n+1),
\end{equation}
and
\begin{equation}
\label{eqMSFt}
0=a_0<a_1<\dots<a_n={\frac{1}{2}}.
\end{equation}

\begin{theorem}
\label{Tpol}
Let ${\cal P}$ be an MSF polygonal as above. For
$j=1,\dots,n$, set
$$ I_j=[a_{j-1},a_j]+m_j. $$
If $K^+=I_1\cup\dots\cup I_n$, then $K=K^+\cup K^-$ is a symmetric
wavelet set of $L^2(\R)$, and it is the disjoint union of $2n$ intervals.

Call $K({\cal P})$ the wavelet set associated to $\cal P$. If
${\cal P}_1\ne{\cal P}_2$ are different polygonals, then
 $K({\cal P}_1)\ne K({\cal P}_2)$.
\end{theorem}

\begin{remark}{\bf (1)}
 Condition~(\ref{eqMSFo}) can be expressed
in the following way. If $m_1$ has the decomposition
\[
 m_1=2^s(2t+1),\ s,t\in{\mathbb N}_0,
\]
 then
\[
 \lambda_n=-s-2,\ m_n=t.
\]
   In fact, by~(\ref{eqMSFo}),
$2^s(2t+1)=m_1=2^{-\lambda_n-2}(2m_n+1)$. In particular, there
exists a bijection between the values of $m_1$ and the couples
$(P_1,P_n)$ that verify~(\ref{eqMSFo}).

\noindent{\bf (2)} Geometrically, condition~(\ref{eqMSFt}) says
that the straight lines $[P_1,P_2],[P_2,P_3]$,
$\dots,[P_{n-1},P_n]$ must have negative, decreasing slopes in
$(-\frac{1}{2},0)$.
\end{remark}
\begin{example}
\label{example}
 Consider, for instance, $n=3$ and nonnegative
integers $s,t,v$ such that $t\ge1$ and $2^v>(2t+1)2^{s+2}$. Let
\[
 P_1=P[0,2^s(2t+1)]=\left(1,2^s(2t+1)\right),\
 P_2=P[-v,0]=\left(2^v,0\right),
\] and
\[
 P_3=P[-s-2,t]=\left(2^{s+2},t2^{s+2}\right).
\]
 Then,~(\ref{eqMSFt}) is satisfied with
\[
 a_1=\frac{2^s(2t+1)}{2^v-1},\ a_2=\frac{2^{s+2}t}{2^v-2^{s+2}}.
\]
 Also~(\ref{eqMSFo}) is satisfied, since
 $\lambda_1=0,\ \lambda_3=-s-2,\ m_1=2^s(2t+1)$ and $m_3=t$.
%$m_1=2^s(2t+1)$, $m_3=t$,
%and
%$$ m_1=2^s(2t+1),\ m_2=0,\ m_3=t. $$
We have, then, the wavelet set
\begin{eqnarray*}
K = K(s,t,v)
 & = &
 \pm\Bigl[{2^s(2t+1)},\frac{2^{s+v}(2t+1)}{2^v-1}\Bigr]\cup
 \Bigl[\frac{2^s(2t+1)}{2^v-1},\frac{2^{s+2}t}{2^v-2^{s+2}}\Bigr] \\
 &   & \cup\Bigl[\frac{2^vt}{2^v-2^{s+2}},\frac{2t+1}{2}\Bigr].
\end{eqnarray*}
%\begin{figure}
%\label{msfone}
%\begin{center}
%\psfrag{A}{$P_1$}
%\includegraphics[height=7cm]{msf2.eps}
%\caption{the MSF polygonal $P_1=(1,1),\ P_2=(\frac{1}{8},\frac{9}{8}),\ P_3=(4,0)$}
%\end{center}
%\end{figure}
\end{example}

\proof[Proof of Theorem~{\rm\ref{Tpol}}]
 It is clear that, if~(\ref{eqMSFt}) holds, then
$\coprod_{m\in{\mathbb Z}}(K^+ +m)$ coincides $a.e.$ with the set
of the reals which are congruent to a number in $[0,\frac{1}{2}]$,
modulo ${\mathbb Z}$. Condition (T) follows by symmetry. Let now
$H_j=2^{-\lambda_j}I_j$, $j=1,\dots,n$. Then, $H_j$ is adjacent to
$H_{j+1}$, for $j=1,\dots,n-1$ since, by definition of $a_j$,
 $$ 2^{-\lambda_j}(a_{j}+m_j)=2^{-\lambda_{j+1}}(a_{j}+m_{j+1}). $$
Moreover, the right endpoint of $H_n$ is twice the left endpoint
of $H_1$, by~(\ref{eqMSFo}). Hence,
 $$ \coprod_{m\in{\mathbb Z}}2^mK^+=[0,\infty)\ a.e. $$
As above, (D) follows by symmetry.
\qed

%The essential features of the wavelet set $K({\cal P})$ can be
%read directly from the picture of ${\cal P}$. Let $P_j=(x_j,y_j)$.
%The position of $I_j\subset[m_j,m_j+\frac{1}{2}]$ is
%approximatively given by $m_j=\frac{y_j}{x_j}$, the slope of the
%straight line through $P_j$ and the origin.
%%==============================================================
\section{Classifying symmetric wavelet sets of $L^2(\R)$}

In this section, we propose a classification of the symmetric
wavelet sets of $L^2(\R)$, which is the basis of the construction 
given in section~2. To each $n$, to each $n$-tuple
$\epsilon=(\epsilon_1,\dots,\epsilon_n)\in\{-1,1\}^n$, and to each
permutation $\tau\in\Sigma_{n-1}$ on the set $\{1,2,\dots,n-1\}$,
we associate a family ${\cal M}(n,\epsilon,\tau)$ of symmetric
wavelet sets of $L^2(\R)$. The main problem of the classification is to
understand what families, among all, are not empty. A partial
result is contained in Proposition~\ref{propminus}.

\begin{proposition}
\label{lemmamain}
 Let $n\in {\mathbb N}$ and let $K\subset{\mathbb R}$
 be a symmetric set, with $K^+$ as in~{\rm(\ref{eqKplus})}.
Then, $K$ is a wavelet set of $L^2(\R)$ if and only if there
exist:
 (i) an $n$-tuple
 $\epsilon=(\epsilon_1,\dots,\epsilon_n)\in\{-1,1\}^n$;
 (ii) a permutation $\tau\in\Sigma_{n-1}$ on the set $\{1,\dots,n-1\}$;
 (iii) a vector $a=(a_0,a_1,\dots,a_n)$ of real numbers that
 satisfy~{\rm(\ref{eqMSFt})};
 (iv) a vector  $m=(m_1,\dots,m_n)$ of nonnegative integers;
 (v) a vector
 $\lambda=(\lambda_1,\dots,\lambda_n)$ of integers; such that
\begin{enumerate}
\item[(a)] $I_j=\epsilon_j[a_{j-1},a_j]+m_j$, $j=1,\dots,n$;
\item[(b)] the intervals $H_k=[\alpha_k,\beta_k]$, $1\le k\le n$, defined by
$H_{\tau(j)}=2^{-\lambda_j}I_j$, $j=1,\dots,n-1$, and $H_n=I_n$,
are such that $2\alpha_1=\beta_n$ and $\beta_k=\alpha_{k+1}$, if
$1\le k\le n-1$;
\item[(c)] $\epsilon$ and $m$ satisfy the condition
\begin{equation}
\label{equnique} 
\mbox{if\ }\epsilon_j=\epsilon_{j+1},\
\mbox{then\ }m_j\ne m_{j+1}.
\end{equation}
\end{enumerate}
Moreover, the wavelet sets relative to $\epsilon\in\{1,-1\}^n$ are
made up of $2n$ disjoint intervals and different sets of data give
different wavelet sets.
\end{proposition}

\proof 
The argument is similar to that in~\cite{HKLS} and we just
give its sketch. For a symmetric $K$, condition (T) is verified if
and only if there are $a$, $m$ and $\epsilon$ such that
~(\ref{eqMSFt}) and (a) hold. On the other hand, (D) holds if and
only if some dyadic dilates of the $I_j$'s have an $a.e.$ disjoint
union of the form $[\alpha,2\alpha]$, for some $\alpha>0$. This is
equivalent to the existence of $\lambda$ and $\tau$ such that (b)
holds. This proves necessity and sufficiency of (a) and (b).

Let now $K$ be a symmetric wavelet set corresponding to a given
$\epsilon\in\{1,-1\}^n$. Then, $K^+$ is made up of {\it at most}
$n$ disjoint intervals $I_1,\dots,I_n$. If two such intervals
share a common endpoint, then there exist two numbers
$0\le\alpha<\beta\le\frac{1}{2}$ and $j,k\in{\mathbb Z}$ such that
either $\alpha+j=\beta+k$, or $\alpha+j=-\beta+k$; or there exists
$j$ such that $\epsilon_j=\epsilon_{j+1}$ and $m_j=m_{j+1}$. The
first case is impossible, the second was ruled out by (c). Hence,
$K^+$ is made up of exactly $n$ disjoint intervals. The uniqueness
part of the last statement follows from this and the fact that the
wavelet set is determined by $a$, $\epsilon$ and $m$. 
\qed

\begin{remark}
\label{remmsf}
\begin{enumerate}
\item[(1)] The set of data $(n,\epsilon,\tau,a,m,\lambda)$ in
  Proposition~\ref{lemmamain} is highly redundant. In fact, it is from
  this redundancy that we obtain the equations that must be satisfied
  by $a$, $\epsilon$ and $m$.
\item[(2)] Without condition (c), it might happen that two
 different sets of data give the
same wavelet set. For instance, the Shannon set is given by both
$(n,\epsilon,\tau,a,m,\lambda)=(1,(-1),I_0,(0,\frac{1}{2}),(1),(0))$
and $(n,\epsilon,\tau,a,m,\lambda)$ = $(2,(-1,-1),I_{1},
(0,\frac{1}{4},\frac{1}{2}),(1,1),(0,0))$, where $I_0$ is a dummy
symbol standing for ``permutation on the empty set'', and
$I_{n-1}$ is the identity function on the set $\{1,\dots,n-1\}$.
\item[(3)] Let ${\cal M}_n(\epsilon,\tau)$ be the family of wavelet sets
relative to a given choice of $(n,\epsilon,\tau)$. If
$(n,\epsilon,\tau)\ne(n^\prime,\epsilon^\prime,\tau^\prime)$, then
${\cal M}_n(\epsilon,\tau)\ne{\cal
M}_{n^\prime}(\epsilon^\prime,\tau^\prime)$.
\item[(4)] Let ${\bf 1}=(1,\dots,1)\in \{-1,1\}^n$.
 Then ${\cal M}_n({\bf 1},I_{n-1})$ is the family of wavelet sets
 considered in section~2.
\end{enumerate}
\end{remark}

The data give different information about the set $K$. The
parameter $a$ says how the classes of reals mod ${\mathbb Z}$ are
divided into the $2n$ intervals of $K$, $\epsilon$ picks those
that will be represented in $K^+$, and $m$ moves them in
$[0,\infty)$.

\medskip

Let ${\bf -1}=(-1,\dots,-1)\in \{-1,1\}^n$. In view of 
Remark~\ref{remmsf} (4) and Theorem~\ref{Tpol}, we might expect
that the wavelet sets in ${\cal M}_n({\bf -1},I_{n-1})$ are
associated to polygonals in the dyadic plane $D$. This is in fact
true, but the only wavelet set obtained in this way is the Shannon
set.

\begin{proposition}
\label{propminus}
\begin{enumerate}
\item[(1)]${\cal M}_1({\bf -1},I_0)$ contains only the Shannon set.
\item[(2)]If $n\ge2$, then ${\cal M}_n({\bf -1},I_{n-1})$ is empty.
\end{enumerate}
\end{proposition}

\proof
 Let $K$ be a subset of the real line. Then $K\in{\cal M}_n({\bf -1},I_{n-1})$,
 with, $K^+=I_1\cup\dots\cup I_n$, if and only if
there exist numbers $b_j$,
\begin{equation}
\label{eqMSFtt}
 \frac{1}{2}=b_0<b_1<\dots<b_n=1
\end{equation}
and $m_j\in{\mathbb N}$, $j=1,\dots,n$, such that
$I_j=[b_{j-1},b_j]+m_j$ and numbers $\lambda_j\in{\mathbb Z}$ such
that, for $j=1,\dots,n-1$,
 $$ 2^{-\lambda_j}(b_j+m_j)=2^{-\lambda_{j+1}}(b_j+m_{j+1}). $$
%and $2(b_0+m_1)=m_n$.
Then, the $b_j$'s, $j=1,\dots,n-1$, are given, as in~(\ref{eqaj}),
by
\begin{equation}
\label{eqbj}
b_j=-\frac{m_j2^{-\lambda_j}-m_{j+1}2^{-\lambda_{j+1}}}
{2^{-\lambda_j}-2^{-\lambda_{j+1}}}.
\end{equation}
Moreover, we can let $\lambda_1=0$, and we have then the condition
\begin{equation}
\label{eqMSFto} 2^{\lambda_n}(2m_1+1)=m_n+1.
\end{equation}

Let $n\ge1$ and ${\cal P}=(P_1,\dots,P_n)$ be the ordered sequence
of points in $D$ given by
$P_j=(2^{-\lambda_j},2^{-\lambda_j}m_j)$, with $m$ and $\lambda$
as defined above. Then, we say that ${\cal P}$ is an {\it MSF($-$)
polygonal}. Thus, there is a bijection between MSF($-$) polygonals
with $n$ vertices and wavelet sets in
 ${\cal M}_n({\bf -1},I_{n-1})$.

Condition~(\ref{eqMSFto}) can be written as follows. Let
$x_n=2^{-\lambda_n},\ y_n=m_n2^{-\lambda_n}$. Then,
 $$ x_n+y_n=2m_1+1. $$

  If $n=1$, we just have to find ${\cal P}=(P_1)$
such that $$ 2m_1+1=m_1+1 $$ which corresponds to the Shannon
wavelet. This proves part (1) of the proposition.

Let $n\ge 2$. We show that there are no MSF($-$) polygonals.
Choose $P_1$ and $P_n$ such that~(\ref{eqMSFto}) holds. Let
$O=(0,0)$, $A=(2m_1+1,0)$ and $B=(0,2m_1+1)$. Then, $P_1=(1,m_1)$
is an interior point of the triangle $AOB$ and $P_n$ belongs to
$AB$, the straight line through $A$ and $B$, by~(\ref{eqMSFto}).

For each couple $P,Q$ of distinct points, let $s(PQ)$ be the
absolute value of $PQ$'s slope. Hence
$s(P_jP_{j+1})\in(\frac{1}{2},1)$. Let $<$ be the ordering on the
line $AB$ for which $B<A$. Now, $s(P_1,A)=\frac{1}{2}$, hence,
$P_2$ still belongs to $AOB$ and $P_1P_2$ intersects $AB$ in a
point $Q_1>A$ having negative ordinate. Since
$\frac{1}{2}<s(P_1,P_2)<s(P_2,P_3)<1$, $P_2P_3$ intersects $AB$ in
a point $Q_2>Q_1$. Inductively, if $P_jP_{j+1}$ intersects $AB$ in
$Q_j$, we have $Q_1<Q_2<\dots<Q_j$, hence all $Q_j$'s have
negative ordinates. But, $P_{n-1}P_n$ intersects the line $AB$ at
$P_n$ which has positive ordinate, since $n\geq 2$. This
contradiction shows that there are no MSF($-$) polygonals.
\qed
%%===============================================================
\section{A geometric construction of $H^2$-wavelet sets}

In this section, we extend the geometric costruction of section~2 for
constructing some $H^2$-wavelet sets. We again consider the set $D$, in the first quadrant of the Cartesian plane, of points $P$ such that
$$
P \equiv
P[\lambda,m]=(2^{-\lambda},2^{-\lambda}m), \quad {\rm where}
~m\in\N_0 ~{\rm and}~ \lambda\in\Z.
$$
 Let $n\in\N$ and
$P_j=P[\lambda_j,m_j]$, $j=1,2,\dots,n$
($m_j\not=m_{j+1},m_0\not=m_n+1$), and
$\tilde{P_n}=P[\lambda_n+1,m_n+1]$. Without loss of generality, we
can take $\lambda_1=0$ and $m_1=0$. Observe that $\tilde{P_n}$ is
uniquely determined once $P_n$ is given. Define the points $a_j$,
$j=1,\dots,n$ as follows:

\begin{equation}
\label{h2-eqa0} a_0=-\frac{(m_n+1)2^{-(\lambda_n+1)}}
{2^{-(\lambda_n+1)}-1},
\end{equation}

\begin{equation}
\label{h2-eqaj}
a_j=-\frac{m_j2^{-\lambda_j}-m_{j+1}2^{-\lambda_{j+1}}}
{2^{-\lambda_j}-2^{-\lambda_{j+1}}}, \quad j=1,2,\dots,n-1,
\end{equation}
and
\begin{equation}
\label{h2-eqan} a_n=a_0+1.
\end{equation}
That is, $a_0$ is the negative of the slope of the straight line
joining $P_1$ and $\tilde{P_n}$, and $a_j$ $(j=1,2,\dots,n-1)$ is
the negative of the slope of the straight line through $P_j$ and
$P_{j+1}$. The ordered sequence of points ${\mathcal
P}=(P_1,\dots,P_n)$ is said to be an {\it $H^2$-MSF polygonal} if
the points $a_j$ $(j=0,1,\dots,n)$ satisfy
\begin{equation}
\label{h2-eqMSFo} 0<a_0<1
\end{equation}
and
\begin{equation}
\label{h2-eqMSFt} a_0<a_1<\dots<a_n=a_0+1.
\end{equation}

\begin{theorem}
\label{h2-Tpol} Let ${\cal P}$ be an $H^2$-MSF polygonal as above.
Let $$ I_j=[a_{j-1},a_j]+m_j,\quad j=1,2,\dots,n. $$ Then
$K=I_1\cup\dots\cup I_n$ is an $H^2$-wavelet set and it is the
disjoint union of $n$ intervals.

 Denote the $H^2$-wavelet set associated to
$\cal P$ by $K({\cal P})$. If ${\cal P}_1\ne{\cal P}_2$ are different polygonals,
then $K({\cal P}_1)\ne K({\cal P}_2)$.
\end{theorem}

The proof of the above theorem is similar to the $L^2$ case (see
Theorem~\ref{Tpol}) and can be omitted.

\begin{example}
Let $n=2$. Consider $\lambda_1=0,m_1=0,\lambda_2=r$ and $m_2=k$,
where $r$ and $k$ are integers. Then we have
\[
P_1=(1,0),\ P_2=(2^{-r},2^{-r}k) ~{\rm and}~
\tilde{P_2}=\bigl(2^{-(r+1)},2^{-(r+1)}(k+1)\bigr).
\]
This gives us
\[
a_0=\frac{k+1}{2^{r+1}-1},\ a_1=\frac{k}{2^r-1}{\rm~and~}
a_2=a_0+1=\frac{k+2^{r+1}}{2^{r+1}-1}.
\]
In order to prove that $(P_1,P_2)$ is an $H^2$-MSF polygonal, we have to
verify~(\ref{h2-eqMSFo}) and~(\ref{h2-eqMSFt}).
Inequality~(\ref{h2-eqMSFo}) is equivalent to $0\leq k<2(2^r-1)$
whereas~(\ref{h2-eqMSFt}) is trivially satisfied. The condition
$m_2\not=m_1$ then gives us $1\leq k < 2(2^r-1)$. Now
\[
I_1=[a_0,a_1]+m_1 = \Bigl[\frac{k+1}{2^{r+1}-1},\frac{k}{2^r-1}\Bigr]
\]
and
\[
 I_2=[a_1,a_2]+m_2 = \Bigl[\frac{2^r k}{2^r-1},
 \frac{2^{r+1}(k+1)}{2^{r+1}-1}\Bigr].
\]
Then
\[
K= \Bigl[\frac{k+1}{2^{r+1}-1},\frac{k}{2^r-1}\Bigr]\cup
\Bigl[\frac{2^r k}{2^r-1},\frac{2^{r+1}(k+1)}{2^{r+1}-1}\Bigr]
\]
is an $H^2$-wavelet set, where $r\geq 1$, $1\leq k < 2(2^r-1)$.
These examples, in fact, exhaust all $H^2$-wavelet sets consisting
of two disjoint intervals (see introduction).
\end{example}
%%=================================================================

\section{Interval wavelet sets of $H^2(\R)$}
The classification of interval wavelet sets of $H^2(\R)$,
analogous to that in Proposition~\ref{lemmamain}, is much simpler.
In this case we prove the following.

\begin{proposition}\label{thm:hardy.int}
Let $I_r=[p_r,q_r], 0\leq r\leq m$ and $K = I_0\cup
I_1\cup\cdots\cup I_m$, with $0<p_0<q_0<p_1<q_1<\cdots <p_m<q_m$.
Then $ K $ is a wavelet set for $H^2(\R)$ if and only if
\begin{enumerate}
\item[(i)] $I_0\cup(I_{\rho(1)}-k_1)\cup(I_{\rho(2)}-k_2)\cup\cdots\cup
(I_{\rho(m)}-k_m) = [p_o,p_0 +1]$, for some permutation $\rho\in
{\Sigma_m},k_i\in\N_0, 1\leq i\leq m$ such that the right endpoint of each interval is equal to the left endpoint of the next
interval; and
\item[(ii)] $I_0\cup(2^{-r_1}I_{\sigma(1)})\cup(2^{-r_2}I_{\sigma(2)})
\cup\cdots\cup(2^{-r_m}I_{\sigma(m)}) = [p_0,2p_0]$, for some
permutation $\sigma\in \Sigma_m, r_i\in\N_0, 1\leq i\leq m$ such
that the right endpoint of each interval is equal to the left
endpoint of the next interval.
\end{enumerate}
\end{proposition}

\proof 
Suppose $K$ is as in the hypothesis and $(i)$ and $(ii)$
hold. Then,
\begin{eqnarray*}
 \bigcup_{k\in\Z}(K+k)
 & = & \bigcup_{k\in\Z}\bigcup_{i=0}^m (I_i+k) \\
 & = & \bigcup_{k\in\Z}\Bigl[\bigl\{I_0\cup\bigcup_{i=1}^m
       (I_{\rho(i)}-k_i)\bigr\}+k\Bigr] \\
 & = & \bigcup_{k\in\Z}\bigl\{[p_0,p_0+1]+k\bigr\} \\
 & = & \R.
\end{eqnarray*}
Also, $\{K+k:k\in\Z\}$ is a pairwise disjoint collection.
Similarly, $\cup_{j\in\Z}(2^j K) = \R^+$ and $\{2^j K:j\in\Z\}$ is
pairwise disjoint. Hence, $K$ is a wavelet set for $H^2(\R)$, by
Theorem~\ref{thm:h2wav}.

 Conversely, suppose $K = I_0\cup\cdots\cup I_m$ is an $H^2$-wavelet set,
where $I_r = [p_r,q_r], 0\leq r \leq m$. So, (T') and (D') of
Theorem~\ref{thm:h2wav} are satisfied. Note that for each
 $r,1\leq r \leq m,$ there is a unique $k_r\in\Z$ such that
 $p_r-k_r\in [p_0,p_0+1]$. But, $p_r-k_r\not\in [p_0,q_0]$.
Otherwise, $I_0\cap (I_r-k_r)$ will contain  the
interval $[p_r-k_r, q_0]$ if $q_0\leq q_r-k_r$, or the interval
$[p_r-k_r, q_r-k_r]$ if $q_r-k_r\leq q_0$. This will be a
contradiction to (T'). So, $ p_r-k_r \in [q_0, p_0+1]$. Note that
$k_r \geq 0$ as $p_r > q_o$.

\smallskip
\noindent{\bf Claim.}
 $I_r-k_r = [p_r,q_r]-k_r\subset [q_0,p_0+1]$
\smallskip

 If not, then $q_0\leq p_r-k_r < p_0+1 < q_r-k_r$. That is,
 $p_r-k_r-1< p_0 < q_r-k_r-1$. But, this implies that
 $I_0\cap(I_r-k_r-1)$ contains the interval $[p_0, q_r-k_r-1]$
 if $q_r-k_r-1 \leq q_0$, or the interval $[p_0,q_0]$ if
 $q_0\leq q_r-k_r-1$. In either case, (T') is violated.
 So, the claim is proved.

Now,
\begin{equation}\label{eqn:hardy1}
\begin{split}
 \bigcup_{k\in\Z}(K+k)
 & = \bigcup_{k\in\Z} \bigcup_{r=0}^m(I_r+k) \\
 & = \bigcup_{k\in\Z} \Bigl\{ \bigl(I_0 \cup (I_1-k_1)
       \cup\cdots\cup (I_m-k_m )\bigr)+k \Bigr\}.
\end{split}
\end{equation}

If there is no $ r, 1\leq r \leq m $ such that $ q_0 = p_r-k_r $,
then $(I_1-k_1)\cup\cdots\cup(I_m-k_m )$ will be properly
contained in $[q_0,p_0+1]$. In fact, the set $[p_0,p_0+1]\setminus
\{I_0 \cup (I_1-k_1)\cup\cdots\cup (I_m-k_m )\}$ will have
positive measure which, in turn, will show that~(\ref{eqn:hardy1})
cannot be an $a.e.$ partition of $\R$. So, there is an index $ r,
1\leq r\leq m $ such that
 $q_0 =p_r-k_r$. Further, such an index is unique. For, if there exist
$r,s$ with $1\leq r,s\leq m$ such that $q_0 = p_r-k_r = p_s-k_s$, then
$(I_r-k_r)\cap (I_s-k_s)$ will contain an interval which will
again contradict (T').

Hence, there is exactly one index $i_1$ such that
 $q_0=p_{i_1}-k_1, 1\leq i_1\leq m$ and $ k_1\in\N_0$.

Now,
 $I_0\cup(I_{i_1}-k_1) = [p_0,q_1-k_1]\subset[p_0,p_0+1]$.
 Arguing as above (the role of $ q_0 $ is now taken by
 $q_1-k_1)$, there is exactly one index $i_2$ such that
\begin{align*}
 q_{i_1}-k_1
 & = p_{i_2}-k_2,\quad 1\leq i_2\leq m,~i_2\not=i_1,
 ~k_2\in\N_0. \\
 \intertext{Similarly,}
 q_{i_2}-k_2 & = p_{i_3}-k_3,\quad 1\leq i_3\leq m,
 ~i_3\not=i_1,i_2,~k_3\in\N_0. \\
 \vdots \\
  q_{i_{m-1}}-k_{m-1}
 & = p_{i_m}-k_m,\quad 1\leq i_m\leq m,
 ~i_m\not=i_1,\cdots,i_{m-1},~k_m\in\N_0.
\end{align*}

Now, $ q_{i_m}-k_m $ has to coincide with $p_0+1 $. Otherwise
$[p_0,p_0+1 ]\setminus \{I_0 \cup (I_{i_1}-k_1)\cup\cdots\cup
(I_{i_m}-k_m )\}$ will have positive measure which will contradict
(T'). So, we have proved (i) of the theorem. By considering dilations by powers of 2 of the intervals $ I_r $ and making use of the partition (D')
of $ \R^+ $, we can prove (ii) in a similar manner.
\qed
%%================================================================

\section{Applications}
In this section, we discuss some applications of the methods
outlined in the previous sections.

\subsection{$n=2,3$}
We give an alternative proof of a theorem  in~\cite{HKLS} and a
subfamily of ${\cal M}_3({\bf 1},I_2)$ that, together with the one in
Example~\ref{example}, exhaust ${\cal M}_3({\bf 1},I_2)$.

\begin{proposition}
\label{PropWTW} 
A symmetric wavelet set $K$ has its positive part
consisting of two disjoint intervals if and only if $$
K^+=\left[2^l,\frac{2^{2l+2}}{2^{l+2}-1}\right]
\cup\left[\frac{2^l}{2^{l+2}-1},\frac{1}{2}\right] \label{eqSol}
$$
for some integer $l\geq 0$.
 In particular, the set of the symmetric wavelet sets with $n=2$
 coincides with ${\cal M}_2({\bf 1},I_1)$.
\end{proposition}

\proof 
Observe that, for $n=2$ as for $n=1$, the permutations play
no role. We only have to consider four cases for $\epsilon$. When
$\epsilon=(-1,-1)$, we have no symmetric wavelet set, by
Proposition~\ref{propminus}. For $\epsilon=(1,-1)$ and
$\epsilon=(-1,1)$, elementary arguments show that there are no
symmetric wavelet sets, and so we are left with ${\cal M}_2({\bf
1},I_1)$. We have to look for MSF polygonals ${\cal P}=(P_1,P_2)$.
But condition~(\ref{eqMSFo}) completely determines $P_2$, if $P_1$
is given. If $P_1=(1,2^l(2t+1))$, then in order for~(\ref{eqMSFt}) to
hold, we must have $t=0$, hence
 $$ P_1=(1,2^l),\ P_2=(2^{l+2},0).$$
Hence, $$ a_1=\frac{2^l}{2^{l+2}-1},\ m_1=2^l,\ m_2=0. $$ From
Theorem~\ref{Tpol}, we obtain the family in the statement of the
proposition. 
\qed
\medskip

We have seen in section~2 a three parameter family of symmetric wavelet sets
which is contained in ${\cal M}_3({\bf 1},I_2)$. Another three
parameter  family can be
constructed as follows. Let $s,u,v$ be integers, $s\ge0$, $u>0$ and
$v>0$, such that
\begin{equation}
\label{eqthree}
2^{s+u}<v<\left(2^{s+u}-1\right)\frac{2^s}{2^s-1}.
\end{equation}
Let
$$
P_1=\left(1,2^s\right),\ P_2=\left(2^{-u},v2^{-u}\right),\
P_3=\left(2^{s+2},0\right).
$$
Then, $(P_1,P_2,P_3)$ is an MSF polygonal. The corresponding wavelet set
is
$$
K=K(s,u,v)=\pm\left[{2^s},\frac{v-2^s}{2^u-1}\right]\cup
\left[\frac{2^uv-2^{s+u}}{2^u-1},\frac{v2^{s+u+2}}{2^{s+u+2}-1}\right]
\cup\left[\frac{v}{2^{s+u+2}-1},\frac{1}{2}\right].
$$
Condition~(\ref{eqthree}) ensures that~(\ref{eqMSFt}) holds, 
while~(\ref{eqMSFo}) holds by our choice of $P_1$ and $P_3$.
It is easy to verify that this family and the one in section~2 exhaust the wavelet
sets in  ${\cal M}_3({\bf 1},I_2)$.
%%-------------------------------------------------------------

\subsection {3-interval $H^2$-wavelet sets}
\label{subsec:3int}
Proposition \ref{thm:hardy.int} allows us to give a complete list
of wavelet sets of $H^2(\R)$ which are union of three disjoint
intervals. We shall see in the next subsection that there are in
fact uncountably many $H^2$-wavelet sets.
%To prove this we will
%construct a family of $H^2$-wavelet sets consisting of four
%intervals whose endpoints depend continuously on areal parameter.

In view of Proposition \ref{thm:hardy.int}, we have the following
result for a 3-interval $H^2$-wavelet set.

\begin{corollary}
\label{cor:3int} Let $K=[p_1,q_1]\cup[p_2,q_2]\cup[p_3,q_3]$ such
that $0 < p_1 < q_1 < p_2 < q_2 < p_3 < q_3$. Then, K is a wavelet
set for $ H^2(\R) $ if and only if for some non-negative integers
$r,s,k,l$
\[
\begin{array}{rl}
(1) ~either~ [{\rm T1}]: & [p_1,q_1]\cup ([p_2,q_2]-k )\cup([p_3,q_3]-l )=[p_1,p_1+1] \\
& with~ q_1=p_2-k, q_2-k=p_3-l,q_3-l=p_1+1 \\
or~[{\rm T2}]: & [p_1,q_1]\cup ([p_3,q_3]- l)\cup([p_2,q_2]- k)=[p_1,p_1+1] \\
& with~ q_1=p_3-l, q_3-l=p_2-k,q_2-k=p_1+1, \\
and~
(2)~either~[{\rm D1}]: &[p_1,q_1]\cup 2^{-r}[p_2,q_2]\cup 2^{-s}[p_3,q_3]=[p_1,2p_1] \\
& with~ q_1=2^{-r}p_2, 2^{-r}q_2=2^{-s}p_3, 2^{-s}q_3=2p_1 \\
or~[{\rm D2}]: & [p_1,q_1]\cup 2^{-s}[p_3,q_3]\cup 2^{-r}[p_2,q_2]=[p_1,2p_1] \\
& with~ q_1=2^{-s}p_3, 2^{-s}q_3=2^{-r}p_2, 2^{-r}q_2=2p_1.
\end{array} \]
\end{corollary}

Thus, in order to characterize all $H^2$-wavelet sets consisting
of three disjoint intervals, we have to consider each of the four
cases (T$i$,D$j$), $i,j=1,2$, and determine the values of the
non-negative integers $r,s,k,l$ such that the corresponding
relations (T$i$,D$j$) hold.
%%--------------------------------------
\subsubsection*{THE CASE (T1,D1)}
We have $K = [p_1,q_1]\cup [p_2,q_2]\cup [p_3,q_3]$ and
\begin{equation}\label{eqn:t1d1}
\begin{split}
 q_1 & =  p_2-k, q_2-k=p_3-l,q_3-l=p_1+1, \\
 q_1 & =  2^{-r}p_2,2^{-r}q_2=2^{-s}p_3, 2^{-s}q_3=2p_1.
\end{split}
\end{equation}

Since $q_1 < p_2$, it is necessary that $k\geq 1$. Now,
$q_2-k=p_3-l \Rightarrow p_3-q_2=l-k$. Since $p_3>q_2$, we have
$l>k$. Thus, $l>k\geq 1$. Similarly, $s>r\geq 1.$ Solving the
equations~(\ref{eqn:t1d1}) for $p_i$'s and $q_i$'s, we get,
\[
\begin{array}{ll}
 p_1 = \frac{l+1}{2^{s+1}-1},          & q_1 = \frac{k}{2^r-1}, \\
 p_2 = \frac{2^r k}{2^r-1},            & q_2 = \frac{l-k}{2^{s-r}-1}, \\
 p_3 = \frac{2^{s-r}(l-k)}{2^{s-r}-1}, & q_3 = \frac{2^{s+1}(l+1)}{2^{s+1}-1}.
\end{array}
\]

We have to ensure that $ 0 < p_1 < q_1 < p_2 < q_2 < p_3 < q_3$.
Clearly, $ p_1 > 0, q_1 < p_2 $ and $ q_2<p_3$. The conditions
$p_1<q_1, p_2<q_2$ and $p_3<q_3$ are equivalent to the following
inequalities.

\begin{enumerate}
\item[(i)] $(2^r-1)l < (2^{s+1}-1)k - (2^r-1)$,
\item[(ii)] $(2^r-1)l > (2^s-1)k $,
\item[(iii)] $ (2^{r+1}-1)l < (2^{s+1}-1)k+2(2^s-2^r) $.
\end{enumerate}

It is easy to see that (iii) $\Rightarrow$ (i). So, we have to
consider only (ii) and (iii). Eliminating  $l$ from (ii) and
(iii), we get
 $0\leq k < 2(2^r-1)$.

Thus, to get all the wavelet sets in this case, we proceed as
follows: Fix an integer $r\geq 1$, then we have to consider only
those integers $k$ such that $1\leq k < 2(2^r-1)$. Consider any
such $k$. By taking $s\geq r+1$, we determine all
integers $l$ satisfying (ii) and (iii). Then any such combination
of $r,k,s$ and $l$ will give rise to an $H^2$-wavelet set.

For example, let $r = 1,$ then $k = 1$. If
 $s=2$, then (ii) and (iii) give us $ l > 3$ and $ 3l < 11$. This
implies $l\geq 4$ and $l\leq 3$. So, there is no integer $l$
satisfying (ii) and (iii). But if we take $s=3$, then (ii) and
(iii) imply $l>7$ and $3l<27$. This gives $l=8$. Hence, $r=1, k=1,
s=3$ and $l=8$ give rise to an $H^2$-wavelet set. The
corresponding wavelet set is
\[
[\tfrac{3}{5},1]\cup [2,\tfrac{7}{3}]\cup
[\tfrac{28}{3},\tfrac{48}{5}].
\]
 If $s=4$, then we get, $l > 15 $ and
 $ 3l < 59 \Rightarrow l = 16,17,18,19$.
 So if we take $r = 1, k= 1, s = 4 $, we get $H^2$-wavelet sets
 for each of the $l$'s; $l= 16,17,18,19.$ The wavelet set
 corresponding to $r=1,k=1,s=4$ and
 $l=16 $ is
\[
[\tfrac{17}{31},1]\cup [2,\tfrac{15}{7}]\cup
[\tfrac{120}{7},\tfrac{544}{31}].
\]

A short table of $(r,k,l,s)$ is given in Table 1.
%%---------------------------------------------------
\subsubsection*{THE CASE (T2,D2)}
We have $ K = [p_1,q_1]\cup [p_2,q_2]\cup [p_3,q_3] $ and
\begin{equation}\label{eqn:t2d2}
\begin{split}
 q_1 & = p_3-l, q_3-l=p_2-k,q_2-k=p_1+1, \\
 q_1 & = 2^{-s}p_3,2^{-s}q_3=2^{-r}p_2, 2^{-r}q_2=2p_1.
\end{split}
\end{equation}

Examining the relations among $p_i$'s and $q_i$'s as in the case
(T1,D1), we have, $ l > k\geq 0$ and $s > r \geq 0.$ Solving the
equations~({\ref{eqn:t2d2}) for $p_i$'s and $q_i$'s , we get
\[
\begin{array}{rclrcl}
p_1 & = & \frac{k+1}{2^{r+1}-1}, & q_1 & = & \frac{l}{2^s-1}, \\
p_2 & = & \frac{l-k}{2^{s-r}-1}, &
          q_2 & = & \frac{2^{r+1}(k+1)}{2^{r+1}-1}, \\
p_3 & = & \frac{2^s l}{2^s -1},&  q_3 & = &
    \frac{2^{s-r}(l-k)}{2^{s-r}-1}.
\end{array} \]

Again, we have to ensure that $0<p_1<q_1<p_2<q_2<p_3<q_3$.
Clearly, $0<p_1$. Since $r+1\leq s$, we get
 $p_1<q_1\Rightarrow q_2 <p_3$. Also, $p_3<q_3\Rightarrow q_1<p_2$.
 So, we have to consider the inequalities $p_1<q_1, p_2<q_2, p_3<q_3$.
These conditions are equivalent to the following inequalities:
\begin{enumerate}
\item[(i)] $ (2^{r+1}-1)l > (2^s-1)(k+1) $,
\item[(ii)] $ (2^{r+1}-1)l < (2^{s+1}-1)k+2(2^s-2^r) $,
\item[(iii)] $ (2^r-1)l > (2^s-1)k $.
\end{enumerate}

If $k = 0$, then (iii) is trivially satisfied and we have to
consider only (i) and (ii), and if $k > 0$, then (iii)
$\Rightarrow $ (i). So, one has to consider (ii) and (iii).
Conditions (ii) and (iii) imply $k < 2(2^r-1)$. If $r=0$ then
 $k<0$, which is not possible.
Thus, to get all $H^2$-wavelet sets in this case, we proceed as
follows:

Fix $r\geq 1$ and consider all $k$'s such that
 $0\leq k <2(2^r-1)$. Take $s>r$. If $k=0$, determine all $l$
 satisfying (i) and (ii); and if $k>0$, then determine all $l$
 which satisfy (ii) and (iii).

For example, $r=1,k=0,s=2$ does not give any wavelet set. But if
we take $r=1,k=0$ and $s=3$, then we get $l=3$. The corresponding
wavelet set is
\[
[\tfrac{1}{3},\tfrac{3}{7}]\cup[1,\tfrac{4}{3}]\cup
[\tfrac{24}{7},4].
\]

Observe that when $ k>0 $, the inequalities to be considered are
same as in the case (T1,D1). So, the table for (T1,D1) also works
for (T2,D2) though we will get different wavelet sets. A short
table for the case $k=0$ is given in Table 2.
%%-----------------------------------------------------
\subsubsection*{THE CASE (T2,D1)}
In this case, we have $ s>r>0$ and $l>k\geq 0$. Solving the
equations [T2] and [D1] of Corollary~\ref{cor:3int} for $p_i$'s
and $q_i$'s, we get
\begin{eqnarray*}
 p_1 & = & \tfrac{1}{2^s}[(2^s-1)k-(2^r-1)l+2^s],\\
 q_1 & = & \tfrac{1}{2^r}[(2^{s+1}-1)k-(2^{r+1}-1)l+2^{s+1}], \\
 p_2 & = & (2^{s+1}-1)k-(2^{r+1}-1)l+2^{s+1}, \\
 q_2 & = & \tfrac{1}{2^s}[(2^{s+1}-1)k-(2^r-1)l+2^{s+1}], \\
 p_3 & = & \tfrac{1}{2^r}[(2^{s+1}-1)k-(2^r-1)l+2^{s+1}], \\
 q_3 & = & 2[(2^s-1)k-(2^r-1)l+2^s].
\end{eqnarray*}
Clearly, $q_1<p_2$ and $ q_2<p_3$. The conditions $0<p_1,\
p_1<q_1,\ p_2<q_2$ and $p_3<q_3$ are equivalent to the following
inequalities:
\begin{enumerate}
\item[(i)] $(2^s-1)k + 2^s > (2^r-1)l$,
\item[(ii)] $[2^s(2^{s+1}-1)-2^r(2^s-1)]k + (2^{2s+1}-2^{r+s}) >
        [2^s(2^{r+1}-1)-2^r(2^r-1)]l$,
\item[(iii)] $(2^s-1)(2^{s+1}-1)k + (2^{2s+1}-2^{s+1}) <
     [2^s(2^{r+1}-1)-(2^r-1)]l$,
\item[(iv)] $[2^{r+1}(2^s-1)-(2^{s+1}-1)]k + 2^{s+1}(2^r-1)
   > (2^r-1)(2^{r+1}-1)l$.
\end{enumerate}

One can show that (ii) $\Rightarrow $ (iv) $\Rightarrow$ (i). So,
we have to consider only (ii) and (iii). Also, as in the previous
cases, eliminating $l$ from (iii) and (iv), we get $k< 2(2^r-1)$.

To get all wavelet sets in this case, we apply a similar
procedure adopted in the case (T1,D1). See Table 3 for some
acceptable values of $(r,k,s,l)$.
%%----------------------------------------------------
\subsubsection*{THE CASE (T1,D2)}
Here we have, $s>r\geq 0$, $l>k\geq 1$. Solving [T1] and [D2] for
$p_i$'s and $q_i$'s, we get
\begin{eqnarray*}
     p_1 & = & \tfrac{1}{2^r}[(2^r-1)l-(2^s-1)k+2^r], \\
     q_1 & = & \tfrac{1}{2^s}[(2^{r+1}-1)l-(2^{s+1}-1)k+2^{r+1}], \\
     p_2 & = & \tfrac{1}{2^s}[(2^{r+1}-1)l-(2^{s}-1)k+2^{r+1}], \\
     q_2 & = & 2[(2^{r}-1)l-(2^{s}-1)k+2^{r}], \\
     p_3 & = & (2^{r+1}-1)l-(2^{s+1}-1)k+2^{r+1}, \\
     q_3 & = & \tfrac{1}{2^r}[(2^{r+1}-1)l-(2^{s}-1)k+2^{r+1}].
\end{eqnarray*}

It is clear that $q_1<p_2$. The inequality $q_2<p_3$ holds if
$p_1<q_1$. Now, the conditions $0<p_1,\ p_1<q_1,\ p_2<q_2$ and
$p_3<q_3$ are equivalent to the following inequalities :
\begin{enumerate}
\item[(i)]   $ (2^r-1)l + 2^r > (2^s-1)k $,
\item[(ii)]  $ [2^s(2^{r}-1)-2^r(2^{r+1}-1)]l + (2^{r+s}-2^{2r+1} <
               [2^s(2^{s}-1)-2^r(2^{s+1}-1)]k $,
\item[(iii)] $ [2^{s+1}(2^r-1)-(2^{r+1}-1)]l + 2^{r+1}(2^s-1) >
               (2^s-1)(2^{s+1}-1)k $,
\item[(iv)]  $ (2^r-1)(2^{r+1}-1)l + 2^{r+1}(2^r-1) <
               [2^r(2^{s+1}-1)-(2^s-1)]k $.
\end{enumerate}

The coefficient of $ l $ in (ii) is negative if and only if
$s=r+1$, in which case (ii) is trivially satisfied. If this
coefficient is non-negative, then it can be shown that (iv)
$\Rightarrow $ (ii). Also, it can be shown that (iii)
$\Rightarrow$ (i). So, we need only (iii) and (iv). Since (iii)
and (iv) imply that $ k < 2(2^r-1)$, the case $r=0$ is ruled out
because $k\geq 1$. Examples of few acceptable values of
$(r,k,s,l)$ are given in Table 4.
%%---------------------------------------------------------

%%----------------------------------------------------
\subsection{There are uncountably many wavelet sets of $L^2(\R)$ and $H^2(\R)$}
When $n=2$, there are countably many symmetric wavelet sets of
$L^2(\R)$. As we have seen in \S6.1, they are disjoint union of
intervals whose endpoints have the form $\frac{p}{2^q-1}$, for
some $p\in{\mathbb Z}$ and $q\in{\mathbb N}$. All examples of
symmetric wavelet sets, even for $n>2$, found in the literature
have endpoints of the above form. In fact all known examples, with
the exception of the Shannon set, belong to the class
 ${\cal M}_n({\bf 1},I_{n-1})$, for some $n$, and then, by the results in
section~2, they have endpoints of this particular form. The class of
symmetric wavelet sets, however, is much richer. In fact, for
$n=3$, we exhibit a family of symmetric wavelet sets which depends
on a real parameter. Hence, there are uncountably many symmetric
wavelet sets and, in particular, their endpoints do not need to be
of the form $\frac{p}{2^q-1}$.
\begin{theorem}
\label{Teugenio} ${\cal M}_3((1,-1,1),I_2)$ contains the family
$\{K_a\colon a\in{\mathbb R}, \frac{1}{3}<a<\frac{1}{2}\}$, where
\begin{equation}
\label{eqEugeni}
K_a^+=\left[a,\tfrac{1}{2}\right]\cup\left[1-a,2a\right]
\cup\left[1,2(1-a)\right].
\end{equation}
Moreover, these wavelet sets are associated to MRA wavelets.
\end{theorem}

In particular, there exist symmetric wavelet sets other than the
Shannon set, whose endpoints are all dyadic rational. For
instance,
\[
K_{3/8}^+=
\left[\frac{3}{8},\frac{4}{8}\right]\cup
\left[\frac{5}{8},\frac{6}{8}\right]\cup
\left[\frac{8}{8},\frac{10}{8}\right].
\]
Another such wavelet set is the one which corresponds to
$a=\frac{7}{16}$. We might think of these sets as ``wavelet sets
on the integers''.

\proof
Subtracting $1$ from the second and third interval in~(\ref{eqEugeni}), 
and leaving the first as it is, we obtain three
intervals whose union is
$$
H=\left[a,\tfrac{1}{2}\right]\cup\left[-a,-(1-2a)\right]
\cup\left[0,1-2a\right].
$$ Now,
$H\coprod(-H)=[-\frac{1}{2},\frac{1}{2}]\ a.e.$, hence (T) holds.
Observe that we have~(\ref{eqMSFt}) with $a_1=1-2a$ and $a_2=a$.

On the other hand, multiplying in~(\ref{eqEugeni}) the first  by
$4$ and the second by $2$, we obtain $a.e.$ disjoint intervals
whose union is $[1,2]$, hence, (D) holds, as well.

A necessary and sufficient condition for a wavelet set $K$ to be
associated to an MRA wavelet is the following~(\cite{FW}, Theorem
3.22). Let $K$ be a wavelet set and $$ K^s=\bigcup_{j\ge
1}(2^{-j}K). $$ By (D), the union is disjoint $a.e.$ and, by (T),
$|K^s|=1$. Then, $K$ is associated to an MRA wavelet if and only
if
\begin{equation}
\label{eqMRA} \coprod_{k\in {\mathbb Z}}(K^s+k)={\mathbb R}\ a.e.
\end{equation}
A direct verification shows that
$$
K_a^s=[a-1,-\tfrac{1}{2}]\cup[-a,a]\cup[\tfrac{1}{2},1-a]\ a.e.
$$

Adding $1$ to the first interval and $-1$ to the third, we obtain
three intervals whose $a.e.$ disjoint union is
$[-\tfrac{1}{2},\tfrac{1}{2}]$, hence~(\ref{eqMRA}) holds.
\qed
\vskip 2mm 

The way the sets in~(\ref{eqEugeni}) were found is the following.
Classes of wavelet sets like ${\cal M}_3((1,-1,1),I_2)$ seem to be
promising places to look for unusual wavelet sets, because $a_1$
and $a_2$ satisfy a system of linear equations without zero
coefficients. The family in~(\ref{eqEugeni}) is one of those for
which the determinant associated to the system vanishes.

\vskip 3mm 
In all the four cases of the 3-interval $H^2$-wavelet
sets (see \S~\ref{subsec:3int}), we observed that the systems of linear equations that
determine the endpoints of the intervals have unique solutions.
Moreover, the endpoints depend upon integer parameters. This fact,
in particular, shows that there are countably many 3-interval
$H^2$-wavelet sets.

We now show that there are uncountably many $H^2$-MSF wavelets by
constructing a family of 4-interval $H^2$-wavelet sets such that
some of the endpoints continuously depend on a real parameter.

\begin{theorem}
For $\frac{1}{2} < c <1$, let $K_c=I_1\cup I_2\cup I_3\cup I_4$,
where
\[
 I_1= [1,2c], I_2= [2c+2,c+3],
 I_3=\bigl[c,\tfrac{c+1}{2}\bigr],~{\rm and}~
 I_4= [\tfrac{c+3}{2},2].
\]
Then, $K_c$ is an $H^2$-wavelet set.
\end{theorem} 

\proof The conditions on $c$ imply that $I_i$'s are nonempty
intervals. Subtracting 1 from the intervals $I_1$ and $I_4$, 3
from $I_2$, and leaving $I_3$ as it is, we get four intervals
whose {\it a.e.} disjoint union is $[0,1]$. That is,
\begin{eqnarray*}
 &   & (I_1-1) \cup (I_2-3) \cup I_3 \cup (I_4-1) \\
 & = & [0,2c-1]\cup [2c-1,c]\cup
[c,\tfrac{c+1}{2}]\cup[\tfrac{c+1}{2},1]=[0,1].
\end{eqnarray*}
Therefore, $\coprod_{k\in\Z}(K_c+k)=\R\ a.e$. Now, multiplying
$I_2$ by $2^{-1}$ and $I_3$ by 2, we get
 intervals whose {\it a.e.} disjoint union is $[1,2]$:
\[ I_1\cup (2^{-1}I_2) \cup (2I_3) \cup I_4 =[1,2], \]
which proves that $\coprod_{j\in\Z} (2^jK_c)=\R^+\ a.e$. By
Theorem~\ref{thm:h2wav}, $K_c$ is an $H^2$-wavelet set. 
\qed
\vskip 2mm 

We recently became aware of the paper~\cite{Maj} by G.
Majchrowska, in which she also arrives at similar conclusions by a
different method. The wavelet sets she obtained are also union of
four intervals.

We end this section with the construction of a family of 5-interval 
$H^2$-wavelet sets where some of the endpoints of the intervals depend 
on two independent real parameters. 

\begin{theorem}
Let $\frac{1}{2}<x<y<1$ and $x+1>2y$. That is, $(x,y)$ is in the interior of the triangle with vertices $(\frac{1}{2}, \frac{3}{4})$, $(\frac{1}{2},1)$ and $(1, 1)$. Then
\begin{equation}\label{E.5int}
K_{x,y}=[x,y]\cup[1,2x]\cup[2y,x+1]\cup[y+1,2]\cup[2x+2,2y+2]
\end{equation}
is an $H^2$-wavelet set.
\end{theorem}

\proof
Let us denote the intervals in the right hand side of~(\ref{E.5int}) by $I_1,I_2,\dots,I_5$. The conditions on $x$ and $y$ ensure that these intervals are non-empty. Observe that the intervals $I_1$, $I_4-1$, $I_2$, $I_5-2$, $I_3$ 
are pairwise disjoint and 
$I_1\cup (I_4-1)\cup I_2\cup(I_5-2)\cup I_3 = [x,x+1]$. Similarly, the intervals  
$I_1$, $2^{-1}I_3$, $2^{-2}I_5$, $2^{-1}I_4$, $I_2$ are pairwise disjoint and 
$I_1\cup(2^{-1}I_3)\cup(2^{-2}I_5)\cup(2^{-1}I_4)\cup I_2 = [x,2x]$. 
Hence, by Theorem~\ref{thm:h2wav}, $K_{x,y}$ is an $H^2$-wavelet set. 
%%===============================================================

\section{Wavelet sets of $L^2(\R)$ accumulating in $0$}
By (D) in Theorem \ref{Tmsf}, a wavelet set $K$ cannot contain a
nondegenerate interval containing $0$. It is natural to ask
whether $0$ can be an accumulation point of $K$. The answer is
indeed yes. Some examples of such wavelet sets are the following:

(1) Madych~\cite{Mad} constructed an example of an MSF wavelet
$\psi$ such that $\hat\psi$ does not vanish in any neighbourhood
of the origin so that $0$ is an accumulation point of the
corresponding wavelet set.

(2) Garrig\'os~\cite{gar}, in his Ph. D. thesis, gave an example
of a wavelet set $K\subset [-2,\frac{1}{2}]$ with the same
property.

(3) In~\cite{BGRW}, the authors constructed wavelet sets
$K_{\epsilon}$ for each $\epsilon$, $0<\epsilon\leq \frac{1}{3}$,
such that
$K_\epsilon\subset[-\tfrac{4}{3},\tfrac{4}{3}+\epsilon]$, and $0$
is an accumulation point of $K_\epsilon$.

In this section we construct some families of wavelet sets of $L^2(\R)$ that
accumulate in $0$. First we construct a symmetric wavelet set
having 0 as an accumulation point. Unfortunately, this set is unbounded.
Next, we construct a family of bounded wavelet sets (but, not
symmetric) with the same property. The final construction is of a
family of bounded symmetric wavelet sets with the origin as an
accumulation point, which provides a positive answer to a question
of Garrig\'os.
%%-----------------------------------------------------------------

\subsection{A symmetric wavelet set that accumulates in $0$}
\label{subsec:dctft0} Here, we exhibit an example of a
{\it symmetric} wavelet set having $0$ as accumulation point.
First, we prove that the example does what it is supposed to, then
we will see how it fits into the geometric scheme developed in section~2
and how this sort of examples can be found.

\begin{proposition}
\label{propBra} Let $K$ be the symmetric subset of ${\mathbb R}$
such that
\begin{equation}
\label{eqBra} K^+=\left(\coprod_{n=0}^\infty
  I_n\right)\cup\left(\coprod_{n=0}^\infty J_n\right)
\end{equation}
where
$$
 I_0=\left[\frac{4}{3},\frac{3}{2}\right],\
 I_n=\left[\frac{3\cdot2^{3n+2}}{2^{2n+3}-1},
  \frac{3\cdot2^{3n+1}}{2^{2n+2}-1}\right],\ n\ge1
$$
   and
$$
 J_0=\left[\frac{1}{5},\frac{1}{3}\right],\
 J_n=\left[\frac{3\cdot2^n}{2^{2n+4}-1},
  \frac{3\cdot2^{n-1}}{2^{2n+3}-1}\right],\ n\ge1.
$$
Then, $K$ is a symmetric wavelet set having $0$ as accumulation point.
\end{proposition}

\proof Being the limit of $J_n$'s endpoints, as $n\to\infty$, $0$
is an accumulation point for $K^+$.

Now, we verify (D). For $n\ge1$, let $\lambda_n=n+1$, and let
$\lambda_0=0$. Similarly, for $n\ge0$, let $\mu_n=-(n+2)$.
Consider the intervals $$ H_n=2^{-\lambda_n}I_n,\
L_n=2^{-\mu_n}J_n,\ n\ge0. $$
Then, $H_0=I_0$ and
$$
H_n=\left[\frac{3\cdot2^{2n+1}}{2^{2n+3}-1},
\frac{3\cdot2^{2n}}{2^{2n+2}-1}\right]
$$
 if $n\ge1$. Also, $L_0=[\frac{4}{5},\frac{4}{3}]$ and, if
$n\ge1$,
$$
L_n=\left[\frac{3\cdot2^{2n+2}}{2^{2n+4}-1},
\frac{3\cdot2^{2n+1}}{2^{2n+3}-1}\right].
$$

Given intervals $A=[a,b]$ and $B=[c,d]$, let us write
$A\rightarrow B$ if $b=c$, that is, if $A$ is adjacent to $B$ and
preceeds $B$ on the real line. Then, the intervals $H_n$ and $L_n$
satisfy the relations
$$
 \dots \rightarrow H_{n+1}\rightarrow
 L_n\rightarrow H_n\dots \rightarrow L_1\rightarrow H_1\rightarrow
 L_0\rightarrow H_0.
$$ Thus,
$$
 \coprod_{n\ge0}\left(L_n\cup H_n\right)=
 \left[\lim_{n\to\infty}\frac{3\cdot2^{2n+2}}{2^{2n+4}-1},\frac{3}{2}\right]
=\left[\frac{3}{4},\frac{3}{2}\right]\ a.e.
$$
 This proves (D), since $\frac{3}{2}=2\cdot\frac{3}{4}$.

In order to verify (T), let $l_0=1$, $l_n=3\cdot2^{n-1}$, if
$n\ge1$, and let $A_n=I_n-l_n$. Then,
$$
A_0=\left[\frac{1}{3},\frac{1}{2}\right],\
A_n=\left[\frac{3\cdot2^{n-1}}{2^{2n+3}-1},\frac{3\cdot2^{n-1}}{2^{2n+2}-1}\right],
\ \mbox{if}\ n\ge1,
$$
and we have the adjacency relations

\begin{equation}
\label{eqadj} \dots \rightarrow A_{n+1}\rightarrow J_n\rightarrow
A_n\dots \rightarrow J_1\rightarrow A_1\rightarrow J_0\rightarrow
A_0.
\end{equation}
Hence, $$ \coprod_{n\ge0}\left(J_n\cup A_n\right)=
\left[\lim_{n\to\infty}\frac{3\cdot2^{n}}{2^{2n+4}-1},\frac{1}{2}\right]
=\left[0,\frac{1}{2}\right]\ a.e. $$ and (T) follows. 
\qed
\vskip 2mm 

To find this sort of examples, it is natural to look for
{\it infinite MSF  polygonals}, corresponding to infinitely many
disjoint intervals in $K$, in which infinitely many points lie on
the axis of the abscissas, so that infinitely many intervals are
contained in $[0,\frac{1}{2}]$. We see how this procedure works in
our case, giving an alternative argument for the proof of
Proposition~\ref{propBra}, which also was the path leading to the
wavelet set in~(\ref{eqBra}).

Let $P_0=(1,1)=P[0,1]$,
$P_n=(x_n,y_n)=(2^{-(n+1)},\frac{3}{4})=P[n+1,3\cdot2^{n-1}]$, if
$n\ge1$, and $Q_n=(u_n,v_n)=(2^{n+2},0)=P[-(n+2),0]$, if $n\ge0$.
With the notation of section~3, let $$ \alpha_n=s(P_n,Q_n),\
\beta_n=s(P_{n+1},Q_n). $$ Also, set $\beta_{-1}=\frac{1}{2}$. We
have the relations $$
0<\dots<\alpha_{n}<\beta_{n-1}<\alpha_{n-1}<\cdots<\beta_0<\alpha_0<\beta_{-1}=
\tfrac{1}{2} $$ that replace~(\ref{eqMSFt}). With the same
notation as above, $$ A_n=[\alpha_n,\beta_{n-1}],\
J_n=[\beta_n,\alpha_n]. $$

As in \S2, consider the intervals
$[\alpha_n,\beta_{n-1}]+\frac{y_n}{x_n}$ and
$[\beta_n,\alpha_n]+\frac{v_n}{u_n}$. These are, respectively, the
intervals $I_n$ and $J_n$ in the statement of the proposition.
Their union is $K^+$, where $K$ is our perspective symmetric
wavelet set. Now, condition (T) is satisfied, by construction.
Moreover, since the intervals
$x_n([\alpha_n,\beta_{n-1}]+\frac{y_n}{x_n})=A_n$ and
$u_n([\beta_n,\alpha_n]+\frac{v_n}{u_n})=J_n$ satisfy the
adjacency relations~(\ref{eqadj}), to verify condition (D), it
suffices to check that $$
2\lim_{n\to\infty}x_n(\alpha_n+\tfrac{y_n}{x_n})=
x_0(\beta_{-1}+\tfrac{y_0}{x_0})=\tfrac{3}{2} $$ which in fact
holds. This last condition replaces~(\ref{eqMSFo}), in the present
example.
%%-----------------------------------------------------------------

\subsection{A family of bounded wavelet sets accumulating in 0}
\label{subsec:dctft}
To construct new wavelet sets from the known
ones, the concepts of dilation equivalence and translation
equivalence of sets will be useful.

\begin{definition}~(\cite{HW})
A measurable set $A$ is said to be {\it translation equivalent} to
a measurable set $B$ if there exists a measurable partition
$\{A_n\}$ of $A$ and $k_n\in\Z$ such that
 $\{A_n+k_n\}$ is a partition of $B$. Similarly,
a measurable set $A$ is {\it dilation equivalent} to a measurable
set $B$ if there exists a measurable partition $\{A_n'\}$ of $A$
and $j_n\in\Z$ such that $\{2^{j_n}A_n'\}$  is a
partition of $B$.
\end{definition}

As a consequence of Theorem~\ref{Tmsf}, we have the following
\begin{corollary}~{\rm (\cite{HW})}
\label{cor:wavset} Let $K$ and $W$ be subsets of $\R$, and $W$ is
both translation and dilation equivalent to $K$. Then $W$ is a
wavelet set if and only if $K$ is so.
\end{corollary}

We now construct a family of bounded wavelet sets having the
origin as an accumulation point so that the associated wavelets
are band-limited and their Fourier transforms are discontinuous at
the origin. Recall that a function is called band-limited if its
Fourier transform has compact support.

Let $n\geq 2.$ Put
\[
\begin{array}{ll}
a_{n} = \frac{2^{n-2}}{2^n-1}, &  b_n = 2a_{n} = \frac{2^{n-1}}{2^n-1}, \\
c_{n} = \frac{2^{n-1}(2^{n-1}-1)}{2^n-1}, & d_{n} = \frac{2^{2n-2}}{2^n-1},\\
e_n=\frac{2^{n-1}-1}{2^n-1}. &
\end{array}
\]
Let 
\[ 
L = [-d_n,-c_n],\ M= [-e_n,-a_n] \quad{\rm and}\quad R =[a_n,b_n]. 
\] 
Observe that $[-e_n,-a_n]\cup\left([-d_n,-c_n]+2^{n-2}\right)\cup[a_n,b_n]
 = [-e_n,b_n]$ so that 
 $\coprod_{k\in\Z}([-e_n,b_n]+k)=\R$. 
 
 Similarly, 
$\coprod_{j\in\Z} 2^j[a_n,b_n]=\R^+$ and 
$\coprod_{j\in\Z}2^j\bigl([-e_n,-a_n]\cup
2^{-n}[-d_n,-c_n]\bigr)=\R^{-}$.
Therefore, $W_n=L\cup M \cup R$ is a wavelet set, 
by Theorem~\ref{Tmsf}.

Let $\epsilon$ be a real number such that $0<\epsilon
<\frac{a_n}{2}$. Define the following sets:
%\begin{eqnarray*}
%P_1 & = &
%\left[\frac{a_n}{2}+\frac{\epsilon}{2^n},\frac{a_n}{2}+\epsilon\right],\\
%P_2 & = & [a_n+2\epsilon,b_n], \\ P_3 & = & [d_n,d_n+2\epsilon].
%\end{eqnarray*}
%\begin{eqnarray*}
\[
P_1  =
\left[\frac{a_n}{2}+\frac{\epsilon}{2^n},\frac{a_n}{2}+\epsilon\right],\
P_2  =  [a_n+2\epsilon,b_n],\ P_3 = [d_n,d_n+2\epsilon].
\]
  To make $P_2$ a nonempty set we need $2\epsilon< a_n$. Let
\[
\begin{array}{ll}
X_0 = P_1-2^{n-2},     & Y_0 = \frac{1}{2^n}X_0, \\ X_l =
Y_{l-1}-2^{n-2}, & Y_l = \frac{1}{2^{n+l}}X_l,\quad l\geq1.
\end{array}
\]

Now define the set
\begin{equation} \label{eqn:Wn}
 W_{n,\epsilon} =
 \Bigl(L\setminus\bigcup_{l=0}^{\infty}X_l\Bigr)\cup
 \Bigl(\bigcup_{l=0}^{\infty}Y_l\Bigr)\cup M
 \cup (P_1\cup P_2 \cup P_3).
\end{equation}

\begin{theorem}
For each $n\geq 2$ and $0<\epsilon< \frac{a_n}{2}$, the set
$W_{n,\epsilon}$ is a bounded wavelet set having $0$ as an
accumulation point.
\end{theorem}

\proof Since $0$ is the limit of the endpoints of the sets $Y_l$,
as $l\to\infty$, it is an accumulation point of $W_{n,\epsilon}$.
We now proceed to prove that $W_{n,\epsilon}$ is a wavelet set. In
view of Corollary~\ref{cor:wavset}, it is enough to show that
$W_{n,\epsilon}$ is translation and dilation equivalent to the
wavelet set $W_n$.

First of all, we show by induction that $X_l\subset L$ for all
$l\geq 0$. Note that $P_1\subset [-a_n,a_n]$. So,
$$
 X_0 = P_1-2^{n-2}\subset[-a_n,a_n]-2^{n-2} =[-d_n,-c_n] =L.
$$ Now, assume that $X_m\subset L$. Then,
\[
Y_m=\frac{1}{2^{m+n}}X_m\subset\left[-\frac{d_n}{2^{m+n}},
-\frac{c_n}{2^{m+n}}\right]\subset[-a_n,a_n].
\]
Therefore,
\[ X_{m+1}=Y_m-2^{n-2}\subset L. \]

The interval $X_0$ lies inside $[-2^{n-2},-c_n]$ and $\{X_l:l\geq 1\}$
lie in $[-d_n,-2^{n-2}]$. Further, $X_{l+1}$ lies to the right of
$X_l,l\geq 1$. The intervals $Y_l,l\geq 0$ lie in
$\frac{1}{2^l}[-a_n,-\frac{a_n}{2}]$. Also observe that
$\{X_l:l\geq 0\}$ and $\{Y_l:l\geq 0\}$ are disjoint collections.

Now, we show that the sets $W_{n,\epsilon}$ and $W_n$ are dilation
equivalent. We have
\begin{eqnarray*}
 2P_1\cup P_2\cup\tfrac{1}{2^n}P_3
 & = & \left[a_n+\frac{\epsilon}{2^{n-1}},a_n+2\epsilon\right]
 \cup[a_n+2\epsilon,b_n]\cup\left[a_n,a_n+\frac{\epsilon}{2^{n-1}}\right] \\
 & = & [a_n,b_n] = R
\end{eqnarray*}
and
\[
 \Bigl(L\setminus\bigcup_{l=0}^{\infty}X_l\Bigr)\cup
 \Bigl(\bigcup_{l=0}^{\infty}2^{n+l}Y_l\Bigr) =
 \Bigl(L\setminus\bigcup_{l=0}^{\infty}X_l\Bigr)\cup
 \Bigl(\bigcup_{l=0}^{\infty}X_l\Bigr)= L.
\]
The last equality follows since $X_l\subset L$ for all
$l\geq 0$. The set $M$ appears in both the partitions of $W_n$ and
$W_{n,\epsilon}$.

To see the translation equivalence of the sets $W_{n,\epsilon}$
and $W_n$, observe that
\[
 P_2\cup(P_3-2^{n-2})=[a_n+2\epsilon,b_n]\cup[a_n,a_n+2\epsilon]=R
\] and
\[
 \Bigl(L\setminus\bigcup_{l=0}^{\infty}X_l\Bigr)\cup
 \Bigl(\bigcup_{l=0}^{\infty}(Y_l-2^{n-1})\Bigr)\cup (P_1-2^{n-2}) =
\]
\[
 \Bigl(L\setminus\bigcup_{l=0}^{\infty}X_l\Bigr)\cup
 \Bigl(\bigcup_{l=1}^{\infty}X_l\Bigr)\cup X_0 =  L.
\]
Again, $M$ appears in both the partitions of
$W_n$ and $W_{n,\epsilon}$. Since the set $W_{n,\epsilon}$ is
translation and dilation equivalent to the wavelet set $W_n$, we
have proved that $W_{n,\epsilon}$ is a wavelet set. \qed
%%-------------------------------------------------------------

\subsection{A family of bounded symmetric wavelet sets accumulating in 0}
\label{subsec:dctfte}
In~\cite{gar}, the question of the existence
of a wavelet $\psi$ of $L^2(\R)$ satisfying the following three
properties was asked:

\begin{itemize}
\item[(i)] $\psi $ is band-limited, i.e., $\hat\psi$ has compact support,
\item[(ii)] $\hat\psi$ is even, and
\item[(iii)] $\hat\psi$ does not vanish in any neighbourhood of the origin.
\end{itemize}

Examples of wavelets satisfying any two of the above three
properties can be constructed. For example, the Shannon wavelet
satisfies (i) and (ii) but not (iii). The wavelet corresponding to 
the wavelet set constructed in \S~\ref{subsec:dctft0} satisfies
(ii) and (iii) but not (i). The wavelets associated with the
wavelet sets of \S~\ref{subsec:dctft}, as well as the wavelet sets
of~\cite{Mad},~\cite{gar} and~\cite{BGRW} referred above, satisfy
(i) and (iii) but not (ii).

In this section we construct a family of wavelets having all the
three properties listed above. These wavelets are again MSF
wavelets and so we construct the associated wavelet sets. To get
these wavelet sets we will suitably modify the method presented in
\S~\ref{subsec:dctft}.

For $n\geq 2$, let $a_n$, $b_n$, $c_n$, $d_n$ and $e_n$ be as in
\S~\ref{subsec:dctft}. Let
\[ L_1 = [-\tfrac{1}{2},-a_n],\ L_2 = [-d_n,-2^{n-2}],\
R_1 = [a_n,\tfrac{1}{2}],\ {\rm and}~R_2 = [2^{n-2},d_n]. \] Let
\[
 K_n = L_1 \cup L_2 \cup R_1 \cup R_2.
\]
  By a simple calculation, we have
\[
 (L_2+2^{n-1})\cup R_2\cup (L_1+(2^{n-2}+1))\cup (R_1+2^{n-2})
 = [c_n,c_n+1].
\]
Hence, $\coprod_{k\in\Z}(K_n+k)=\R$. Also, observe that $(2^n
R_1)\cup R_2 = [\frac{d_n}{2},d_n]$. Hence, $\coprod_{j\in\Z}
2^j(R_1\cup R_2)=\R^+$. By symmetry, $\coprod_{j\in\Z} 2^j(L_1\cup
L_2)=\R^-$, which shows that $\coprod_{j\in\Z} 2^jK_n=\R$.
Therefore, by Theorem~\ref{Tmsf}, $K_n$ is a wavelet set.

For $0<\epsilon<\frac{1}{4}\bigl(\frac{2^{n-1}-1}{2^n-1}\bigr) =
\frac{e_n}{4}$, construct the following sets:
%\begin{eqnarray*}
%S_1 & = &
%\left[\frac{a_n}{2}+\frac{\epsilon}{2^n},\frac{a_n}{2}+\epsilon\right],\\
%S_2 & = & [a_n+2\epsilon,\tfrac{1}{2}],\\ S_3 & = &
%[d_n,d_n+2\epsilon]
%\end{eqnarray*}
\[
S_1 = \left[\frac{a_n}{2}+\frac{\epsilon}{2^n},
\frac{a_n}{2}+\epsilon\right],\ S_2 =
[a_n+2\epsilon,\tfrac{1}{2}],\ S_3 = [d_n,d_n+2\epsilon] \] and
\[
 T_i = -S_i\quad{\rm for}\quad i=1,2,3.
\]
To make $S_2$ a nonempty set, we need to take
$a_n+2\epsilon<\frac{1}{2}$ which is equivalent to
$\epsilon<\frac{e_n}{4}$. Let
\[
\begin{array}{ll}
E_0 = S_1+2^{n-2}, &    F_0 = \frac{1}{2^{n+1}}E_0, \\ E_l =
F_{l-1}+2^{n-2}, & F_l = \frac{1}{2^{n+l+1}}E_l,~~ l\geq1, \\
\end{array}
\]
\[
 G_l = -E_l\quad\mbox{and}\quad H_l = -F_l,\quad\mbox{for all}~l\geq 0.
\]
Define
\begin{eqnarray*}
 K_{n,\epsilon}
 & = & \Bigl(R_2\setminus\bigcup_{l=0}^{\infty}E_l\Bigr)\cup
 \Bigl(\bigcup_{l=0}^{\infty}F_l\Bigr)
 \cup (S_1\cup S_2 \cup S_3) \\
 &   & \cup\Bigl(L_2\setminus\bigcup_{l=0}^{\infty}G_l\Bigr)\cup
 \Bigl(\bigcup_{l=0}^{\infty}H_l\Bigr)\cup (T_1\cup T_2 \cup T_3).
\end{eqnarray*}

\begin{theorem}
For each $n\geq 2$ and $0<\epsilon<\frac{e_n}{4}$, the set
$K_{n,\epsilon}$ is a bounded symmetric wavelet set having $0$ as
an accumulating point.
\end{theorem}

\proof
Clearly, the origin is an accumulation point of the set
$K_{n,\epsilon}$, being the limit (as $l\rightarrow\infty$) of the  endpoints of the
intervals $F_l$, $l\geq 0$. As in the previous theorem, it is
enough to show that $K_{n,\epsilon}$ is translation and dilation
equivalent to the wavelet set $K_n$. Again, we can show  by
induction that $E_l\subset R_2,$ for all $l\geq 0$. By  symmetry
of the set $K_{n,\epsilon}$, it follows that $G_l\subset L_2$ for
all $l\geq 0$.

The intervals $E_l$, $l\geq 0$ lie inside the interval
$[2^{n-2},d_n]$ and $E_{l+1}$ lies to the left of $E_l$ for all
$l\geq 0$. Similarly, the intervals $F_l$, $l\geq 0$ lie in
$\frac{1}{2^l}\bigl[\frac{a_n}{4},\frac{a_n}{2}\bigr]$ so that
$F_{l+1}$ lies to the left of $F_l$ for all $l\geq 0$. Similar
statements are true for the intervals $G_l$ and $F_l$, $l\geq 0$.

We have
\begin{eqnarray*}
 2S_1\cup S_2\cup\tfrac{1}{2^n}S_3
 & = & \left[a_n+\frac{\epsilon}{2^{n-1}},a_n+2\epsilon\right]
 \cup\Bigl[a_n+2\epsilon,\tfrac{1}{2}\Bigr]
 \cup\left[a_n,a_n+\frac{\epsilon}{2^{n-1}}\right] \\
 & = & [a_n,\tfrac{1}{2}] = R_1
\end{eqnarray*}
and
\[
 \Bigl(R_2\setminus\bigcup_{l=0}^{\infty}E_l\Bigr)
 \cup\Bigl(\bigcup_{l=0}^{\infty}2^{n+l+1}F_l\Bigr) =
 \Bigl(R_2\setminus\bigcup_{l=0}^{\infty}E_l\Bigr)\cup
 \Bigl(\bigcup_{l=0}^{\infty}E_l\Bigr) = R_2.
\]
Similarly for $L_1$ and $L_2$. This proves the dilation
equivalence of $K_{n,\epsilon}$ and $K_n$.

The translation equivalence of the sets $K_{n,\epsilon}$ and $K_n$
follows from the following observation. $$
 S_2\cup(S_3-2^{n-2})=
 [a_n+2\epsilon,\tfrac{1}{2}]\cup[a_n,a_n+2\epsilon]=R_1,
$$ and
$$
 \Bigl(R_2\setminus\bigcup_{l=0}^{\infty}E_l\Bigr)
 \cup\Bigl(\bigcup_{l=0}^{\infty}(F_l+2^{n-2})\Bigr)
 \cup(S_1+2^{n-2})
$$
$$ = \Bigl(R_2\setminus\bigcup_{l=0}^{\infty}E_l\Bigr)
\cup\Bigl(\bigcup_{l=1}^{\infty}E_l\Bigr)\cup E_0 = R_2.
$$
Similarly for $L_1$ and $L_2$. Therefore, $K_{n,\epsilon}$ is a
wavelet set.
\qed

\medskip

Let $\hat\psi_{n,\epsilon}$ be the characteristic function of the
set $K_{n,\epsilon}$ . Then, $\psi_{n,\epsilon}$ is a band-limited
wavelet such that $\hat\psi_{n,\epsilon}$ is even and does not
vanish in any neighbourhood of the origin.

\begin{remark}
The presence of the real parameter $\epsilon$ in the family of the
wavelet sets $K_{n,\epsilon}$ proves, in particular, that there
are uncountably many symmetric wavelet sets of $L^2(\R)$, a fact
already proved in Theorem~\ref{Teugenio}.
\end{remark}

\begin{footnotesize}
\begin{center}
%\begin{table}
\begin{tabular}[t]{|r|l|r|c|}
 \multicolumn{4}{c}{Table 1.} \\
 \multicolumn{4}{c}{(T1,D1)}  \\ \hline
$r$ & $k$ & $s$ & $l$        \\ \hline \hline
 1  &  1  &  2  &  --        \\
    &     &  3  &  8         \\
    &     &  4  & 16--19     \\
    &     &  5  & 32--40     \\
    &     &  6  & 64--83     \\
    &     &  7  & 128--168   \\
    &     &  8  & 256--339   \\
    &     &  9  & 512--680   \\
    &     & 10  & 1024--1363 \\ \hline
 2  &   1 &  3  & 3          \\
    &     &  4  & 6,7        \\
    &     &  5  & 11--16     \\
    &     &  6  & 22--35     \\
    &     &  7  & 43--71     \\ \cline{2-4}
    &   2 &  3  & 5          \\
    &     &  4  & 11,12      \\
    &     &  5  & 21--25     \\
    &     &  6  & 43--53     \\ \cline{2-4}
    &   5 &  3,4  &  --      \\
    &     &  5  &  52        \\ \hline
  10&   1 & 11  &   3        \\
    &     & 12  &  5--7      \\
    &     & 20  &  1026--2047\\
    &     & 25  &  32801--65567\\ \hline
  20&    1& 25  &  33--63    \\   \hline
\end{tabular}
\begin{tabular}[t]{|l|r|c|}
 \multicolumn{3}{c}{Table 2.} \\
 \multicolumn{3}{c}{(T2,D2), $k=0$} \\ \hline
$r$      & $s$ & $l$ \\ \hline \hline
 1       &  2  &  --        \\
         &  3  &  3         \\
         &  4  & 6--9       \\
         &  5  & 11--19     \\
         &  6  & 22--41     \\
         &  7  & 43--83     \\ \hline
 2       &  3  & --         \\
         &  4  & 3,4        \\
         &  5  & 5--7       \\
         &  6  & 10--17     \\
         &  7  & 19--35     \\ \hline
 5       &  6  &  --        \\
         &  7  & 3          \\
         &  8  & 5--7       \\
         &  9  & 9--15      \\
         & 10  &  17--31    \\ \hline
\end{tabular}
\begin{tabular}[t]{|r|l|r|c|}
 \multicolumn{4}{c}{Table 3.} \\
 \multicolumn{4}{c}{(T2,D1)} \\ \hline
$r$ & $k$ & $s$ & $l$  \\ \hline \hline 1   & 0   &  2  & --   \\
    &     &  3  &  5   \\
    &     &  4  & --   \\
    &     &  5  & 21   \\
    &     &  6  & --   \\
    &     &  7  & 85   \\
    &     &  8  & --   \\ \cline{2-4}
    &  1  &  2  & 5    \\
    &     &  3  & --   \\
    &     &  4  & --   \\
    &     &  5  & --   \\ \hline
  2 &  0  &  3  & 2    \\
    &     &  4  & --   \\
    &     &  5  & 9    \\  \hline
\end{tabular}
\begin{tabular}[t]{|l|l|r|c|}
 \multicolumn{4}{c}{Table 4.} \\
 \multicolumn{4}{c}{(T1,D2)} \\ \hline
$r$ & $k$ & $s$ & $l$   \\ \hline \hline
 1  & 1   & $s$ & $2^s-2$   \\ \hline
 2  &  1  &  3  & --   \\
    &     &  4  & --   \\
    &     &  5  & --   \\
    &     &  6  & 20   \\ \hline
  3 &  1  &  4  & --   \\
    &     &  5  & --   \\
    &     &  6  & 8    \\
    &     &  7  & --   \\
    &     &  8  & --   \\
    &     &  9  & 72   \\ \hline
\end{tabular}
%\end{table}
\end{center}
\end{footnotesize}

%\pagebreak

\vspace{15 pt}


\begin{thebibliography}{BGRW}
\addcontentsline{toc}{section}{References}

\bibitem[Aus]{aus} P. Auscher.
{\it Solution of two problems on wavelets},
J. Geom. Anal., {\bf 5}, no. 2 (1995), pp. 181-236.

\bibitem[BGRW]{BGRW}
L. Brandolini, G. Garrig\'{o}s, Z. Rzeszotnik, and G. Weiss.
{\it The behaviour at the origin of a class of band-limited wavelets},
Contemporary Mathematics, {\bf 247} (1999), pp. 75-91.

\bibitem[FW]{FW} X. Fang, and X. Wang.
{\it Construction of Minimally Supported Frequency Wavelets},
J. Fourier Anal. and Appl., {\bf 2}, no. 4 (1996), pp. 315-327.

\bibitem[Gar]{gar} G. Garrig\'{o}s.
{\it The characterization of wavelets and related functions
and the connectivity of $\alpha$-localized wavelets on $\R$},
Ph.D. Thesis, Washington University, St.louis (1998).

\bibitem[HKLS]{HKLS}
Y. Ha, H. Kang, J. Lee, and J.K. Seo.
{\it Unimodular Wavelets for $L^2$ and the Hardy space $H^2$},
Michigan Math. J., {\bf 41} (1994), pp. 345-371.

\bibitem[HW]{HW} E. Hern\'andez, G. Weiss.
{\it A First Course on Wavelets}, CRC Press (1996).

\bibitem[HWW]{HWW} E. Hern\'andez, X. Wang, and G. Weiss.
{\it Smoothing Minimally Supported Frequency (MSF) Wavelets: Part II},
J. Fourier Anal. Appl. {\bf 3}, no.1 (1997), pp. 23-41.

\bibitem[Mad]{Mad} W. R. Madych.
{\it Some elementary properties of multiresolution analyses of
$L^2(\R^n)$}, in {\it Wavelets-A Tutorial in Theory and
Applications} (C.K. Chui, Ed.), Academic Press (1992),
pp. 259-294.

\bibitem[Maj]{Maj} G. Majchrowska.
{\it Some new examples of wavelets in the Hardy space $H^2(\R)$,}
Bull. Polish Acad. Sci. Math., {\bf 49}, no. 2 (2001),
pp. 141-149.

%\bibitem[SoW]{SoW} P. M. Soardi, and D. Weiland,
%{\it Single wavelets in $n$-dimensions,}
%J. Fourier Anal. and Appl. {\bf 4}, n.3 (1998), pp. 299-315.

\end{thebibliography}
\end{document}